\pdfoutput=1
\RequirePackage{ifpdf}
\ifpdf 
\documentclass[pdftex]{sigma}
\else
\documentclass{sigma}
\fi

\numberwithin{equation}{section}

\newtheorem{Theorem}{Theorem}[section]
\newtheorem*{Theorem*}{Theorem}
\newtheorem{Corollary}[Theorem]{Corollary}
\newtheorem{Lemma}[Theorem]{Lemma}
\newtheorem{Proposition}[Theorem]{Proposition}

\theoremstyle{definition}
\newtheorem{Definition}[Theorem]{Definition}

\newtheorem{Example}[Theorem]{Example}
\newtheorem{Remark}[Theorem]{Remark}

\usepackage{xcolor}
\usepackage[mathscr]{eucal}
\usepackage{mathtools}
\usepackage{bm}
\usepackage[math]{cellspace}
\usepackage[inline]{enumitem}
\setlist[enumerate,1]{label=(\arabic*),ref=\arabic*,itemsep=1pt}
\setlist[itemize,1]{itemsep=1pt}

\usepackage{tikz}
\usetikzlibrary{positioning}
\usetikzlibrary{cd}
\usetikzlibrary{shapes.geometric}
\usetikzlibrary{decorations.pathreplacing,decorations.markings}
\usetikzlibrary{backgrounds}
\usetikzlibrary{arrows.meta}

\DeclareMathOperator{\diag}{diag}

\DeclareMathOperator{\tr}{tr}
\DeclareMathOperator{\mat}{Mat}

\DeclareMathOperator{\Hom}{Hom}

\DeclareMathOperator{\Trop}{Trop}
\DeclareMathOperator{\SL}{SL}
\newcommand*{\maxzero}[1]{[#1]_{+}}

\tikzset{
mid arrow/.style={postaction={decorate,decoration={
 markings,
 mark=at position .55 with {\arrow[#1]{To[length=1.1mm]}}
}}},
mid 7 arrow/.style={postaction={decorate,decoration={
 markings,
 mark=at position .7 with {\arrow[#1]{To[length=1.1mm]}}
}}},
vertex/.style={circle,inner sep=1.5pt},
e_vertex/.style={ellipse,inner sep=1.5pt},
fan arrow/.style=-{Stealth[length=1.9mm,width=1.7mm]},
help lines/.style={thin,draw=black!20},
quiver/.style={thick, -Stealth}
}

\begin{document}
\allowdisplaybreaks

\newcommand{\arXivNumber}{2301.13239}

\renewcommand{\PaperNumber}{094}

\FirstPageHeading

\ShortArticleName{Periodic $Y$-Systems and Nahm Sums: The Rank 2 Case}

\ArticleName{Periodic $\boldsymbol{Y}$-Systems and Nahm Sums:\\ The Rank 2 Case}

\Author{Yuma MIZUNO}

\AuthorNameForHeading{Y.~Mizuno}

\Address{School of Mathematical Sciences, University College Cork, \\ Western Gateway Building, Western Road, Cork, Ireland}
\Email{\mail{ymizuno@ucc.ie}}
\URLaddress{\url{https://yuma-mizuno.github.io/}}

\ArticleDates{Received June 04, 2025, in final form October 26, 2025; Published online November 03, 2025}

\Abstract{We classify periodic $Y$-systems of rank~2 satisfying the symplectic property. We~find that there are six such $Y$-systems. In all cases, the periodicity follows from the existence of two reddening sequences associated with the time evolution of the $Y$-systems in positive and negative directions, which gives rise to quantum dilogarithm identities associated with Donaldson--Thomas invariants. We also consider $q$-series called the Nahm sums associated with these $Y$-systems. We see that they are included in Zagier's list of rank~2 Nahm sums that are likely to be modular functions. It was recently shown by Wang that they are indeed modular functions.}


\Classification{11F03; 11P84; 13F60}

\section{Introduction}
\subsection{Background}
The $Y$-system is a system of algebraic relations satisfied by coefficients of a cluster algebra,
which has the following form:
\begin{equation}\label{eq:y-system}
 Y_i (u) Y_i (u - r_i) =
 \prod_{j \in I} \prod_{p=1}^{r_i-1} Y_j(u-p)^{\maxzero{n_{ij;p}}} (1 + Y_j(u - p) )^{- n_{ij;p} },
\end{equation}
where $I$ is a finite index set, $Y_i (u)$ for $i \in I$, $u \in \mathbb{Z}$ are commuting variables,
$r_i \in \mathbb{Z}_{\geq 1}$, and~${n_{ij;p} \in \mathbb{Z}}$.
We also use the notation $\maxzero{n} \coloneqq \max(0, n)$.
Such equations are first discovered by Zamolodchikov in the study of thermodynamic Bethe ansatz \cite{Zamo},
prior to the discovery of cluster algebras by Fomin and Zelevinsky \cite{FZ1}.
The most striking feature of Zamolodchikov's \mbox{$Y$-systems}, as well as their generalizations \cite{KunibaNakanishi92, RVT}
defined shortly after the Zamolodchikov's work, is that they are periodic, which was fully proved by applying
the theory of cluster algebras~\cite{FZ_Ysystem, FZ4, IIKKNa, IIKKNb, Keller}.

A systematic treatment of the $Y$-systems in the general setting of cluster algebras, including the $Y$-systems arising from
the thermodynamic Bethe ansatz as special cases, was given by Nakanishi \cite{Nakb}.
This approach was further developed in \cite{mizuno2020difference},
and it was shown that the algebraic relation \eqref{eq:y-system} arises from a cluster algebra if and only if
the data $r_i$, $n_{ij;p}$ have a certain symplectic property.
This allows the ``axiomatic'' study of $Y$-systems without explicitly referring to cluster algebras.
In this general setting, however, the $Y$-system is typically not periodic,
and so the study of periodic $Y$-systems as a generalization of Zamolodchikov's $Y$-systems would be further developed.
In particular, the classification problem for periodic $Y$-systems is a challenging open problem (see the last comments in \cite[Section 3]{Nakb}).

There are several classification results in the literature.
Fomin and Zelevinsky \cite{FZ4} showed that the classification when $r_i = 2$, $n_{ij;p} \leq 0$, and $n_{ii;p} = 0$ for any $i$, $j$, $p$
coincides with the Cartan--Killing classification. Galashin and Pylyavskyy \cite{GalashinPylyavskyy} generalized this result to show that
the classification when $r_i = 2$ and $n_{ii;p} = 0$ for any $i$, $p$ coincides with the classification of ADE bigraphs of Stembridge \cite{Stembridge}.
On the other hand, the situation is more complicated when~${r_i > 2}$ for some $i$, and so far there has been no comprehensive classification results
except when $\lvert I \rvert = 1$ where it is not difficult to give a complete classification
thanks to the work by Fordy and Marsh~\cite{FordyMarsh} (e.g., see \cite[Example 5.6]{mizuno2020difference}).

In this paper, we make a first attempt to give a classification result involving the case $r_i > 2$ for some $i$.
Precisely, we classify the periodic $Y$-systems of the form \eqref{eq:y-system} with $\lvert I \rvert = 2$ satisfying the symplectic property.
We would like to emphasize that we consider general $r_i$, $n_{ij;p}$ in the classification.
The result is given in the next section.

We also discuss the relation to Nahm's conjecture on $q$-series \cite{Nahm, Zagier}
in Section \ref{section:nahm sum}.

\subsection{Main result}\label{section:main relsult}
Let $I$ be a finite set.
We denote by $\mathbb{Y}_0$ the set of pairs $(\bm{r}, \bm{n})$ where $\bm{r} = (r_i)_{i \in I}$ and $\bm{n} = (n_{ij;p})_{i, j \in I, p \in \mathbb{N}}$
are families of integers satisfying $r_i \geq 1$ for any $i$ and
\begin{equation}\label{eq:Y1}
n_{ij;p} = 0\qquad \text{unless} \quad 0 < p < r_i
\end{equation}
for any $i$, $j$, $p$.

\begin{Definition}
 Let $(\bm{r}, \bm{n}) \in \mathbb{Y}_0$.
 Let $\mathbb{P}$ be a semifield, and $(Y_i(u))_{i \in I, u \in \mathbb{Z}}$ be a family of
 elements in $\mathbb{P}$.
 We say that $(Y_i(u))$ \emph{satisfies the $Y$-system} associated with the pair $(\bm{r}, \bm{n})$ if the relation~\eqref{eq:y-system}
 holds for any $i,u$.
 The equation \eqref{eq:y-system} itself is called the \emph{$Y$-system} associated with~$(\bm{r}, \bm{n})$.
 We also say that $(Y_i(u))$ is a \emph{solution of the $Y$-system} if it satisfies the $Y$-system.
\end{Definition}

It is useful to think a pair $(\bm{r}, \bm{n}) \in \mathbb{Y}_0$ as a triple of matrices with polynomial entries
by the map
$(\bm{r}, \bm{n}) \mapsto (N_0(z), N_+(z), N_-(z)) \colon \mathbb{Y}_0 \to (\mat_{I \times I} \mathbb{N} [z])^3$ defined by
\[
 N_0(z) \coloneqq \diag\bigl(1 + z^{r_i}\bigr)_{i \in I}, \qquad
 N_{\pm} (z) \coloneqq \biggr(\sum_{p \in \mathbb{N}} n_{ij;p}^\pm z^p \biggl)_{i, j \in I},
\]
where we set $n_{ij;p}^\pm \coloneqq \maxzero{\pm n_{ij;p}}$.
We also define the map $(\bm{r}, \bm{n}) \mapsto A_\pm (z)\colon \mathbb{Y}_0 \to (\mat_{I \times I} \mathbb{Z} [z])^2$
by~${A_{\pm}(z) \coloneqq N_0(z) - N_{\pm}(z)}$.
Since this map is injective by the condition~\eqref{eq:Y1},
we will identify~$\mathbb{Y}_0$ with the image of this map.
For example, we will use the term ``the $Y$-system associated with~${A_\pm (z) \in \mathbb{Y}_0}$''.

\begin{Definition}
 We say that $A_\pm (z) \in \mathbb{Y}_0$ satisfies the \emph{symplectic property} if
 \begin{align}\label{eq:symplectic property}
 A_+ (z) A_-\bigl(z^{-1}\bigr)^{\mathsf{T}} = A_- (z) A_+ \bigl(z^{-1}\bigr)^{\mathsf{T}},
 \end{align}
 where $\mathsf{T}$ is the transpose of a matrix.
 We denote by $\mathbb{Y}$ the subset of $\mathbb{Y}_0$ consisting of pairs satisfying the symplectic property.
\end{Definition}

The pair $A_\pm(z) \in \mathbb{Y}_0$ satisfies the symplectic property if and only if
the $Y$-system associated with $A_\pm (z) \in \mathbb{Y}_0$ is realized as the exchange relations of coefficients
in a cluster algebra \cite{mizuno2020difference}.
We review this fact in Section \ref{section:Y from cluster}.

\begin{table}[t]
 \centering
 \begin{tabular}{ Sl|Sl|Sl|Sl| Sl }
 \hline
 $A_+(z)$ &
 $A_-(z)$ &
 $h_+$ &
 $h_-$ & \\
 \hline
 $\begin{pmatrix}
 1 + z^2 & -z \\
 -z & 1 + z^2
 \end{pmatrix}$ &
 $\begin{pmatrix}
 1 + z^2 & 0 \\
 0 & 1 + z^2
 \end{pmatrix}$ &
 $3$ & $2$ &
 (1) \\
 $\begin{pmatrix}
 1 + z^2 & -z \\
 - z - z^5 & 1 + z^6
 \end{pmatrix}$ &
 $\begin{pmatrix}
 1 + z^2 & 0 \\
 -z^3 & 1 + z^6
 \end{pmatrix}$ &
 $8$ & $6$ &
 (2) \\
 $\begin{pmatrix}
 1 + z^2 & - z \\
 - z - z^5 - z^9 & 1 + z^{10}
 \end{pmatrix}$ &
 $\begin{pmatrix}
 1 + z^2 & 0 \\
 -z^3 - z^7 & 1 + z^{10}
 \end{pmatrix}$ &
 $18$ & $10$ &
 (3) \\
 $\begin{pmatrix}
 1 + z^2 & - z \\
 - z & 1 + z^2
 \end{pmatrix}$ &
 $\begin{pmatrix}
 1 + z^2 - z & 0 \\
 0 & 1 + z^2 - z
 \end{pmatrix}$ &
 $3$ & $3$ &
 (4) \\
 $\begin{pmatrix}
 1 + z^2 & -z \\
 -z -z^2 & 1 + z^3
 \end{pmatrix}$ &
 $\begin{pmatrix}
 1 + z^2 - z & 0 \\
 0 & 1 + z^3
 \end{pmatrix}$ &
 $5$ & $3$ &
 (5) \\
 $\begin{pmatrix}
 1 + z^2 & - z \\
 - z & 1 + z^2 - z
 \end{pmatrix}$ &
 $\begin{pmatrix}
 1 + z^2 & 0 \\
 0 & 1 + z^2
 \end{pmatrix}$ &
 $5$ & $2$ &
 (6) \\
 \hline
 \end{tabular}
\caption{Finite type classification for $Y$-systems of rank $2$. The numbers $h_{\pm}$ are the length of
reddening sequences in positive and negative directions, respectively.}
\label{table:finie type}
\end{table}

\begin{Definition}\label{def:periodic}
 We say that a solution of a $Y$-system is \emph{periodic} if
 there is a positive integer~${\Omega > 0}$ such that $Y_i (u + \Omega) = Y_i (u)$
 for any $i$, $u$.
\end{Definition}

\begin{Definition}
 We say that a pair $A_\pm (z) \in \mathbb{Y}$ is of \emph{finite type} if
 any solution (in any semifield) of the $Y$-system associated with this pair is periodic.
 In this case, we also say that $Y$-system itself is \emph{periodic}.
\end{Definition}

The purpose of this paper is to classify periodic $Y$-systems of rank $2$.
Before stating the result, we give a few remarks.
We say that $A_\pm(z) \in \mathbb{Y}_I$ is \emph{decomposable} if
it is a direct sum of some $A'_\pm(z) \in \mathbb{Y}_{I'}$ and $A''_\pm(z) \in \mathbb{Y}_{I''}$
with nonempty $I'$ and $I''$.
We say that $A_\pm(z) \in \mathbb{Y}_I$ is \emph{indecomposable} if it is not decomposable.
It is enough to consider indecomposable pairs in the classification.
We also note that $A_\pm (z)$ is of finite type if and only if
$A_\pm(z)^{\mathrm{op}} \coloneqq A_\mp(z)$ is of finite type by the correspondence
between solutions $Y_i(u) \mapsto Y_i(u)^{-1}$.
The main results are summarized as follows.

\begin{Theorem}\label{theorem:main}
 Suppose that $I = \{1, 2\}$.
 \begin{enumerate}\itemsep=0pt
 \item[$(1)$] Any pair $A_\pm (z) \in \mathbb{Y}$ in Table {\rm\ref{table:finie type}} is of finite type.
 \item[$(2)$] Any indecomposable pair $A_\pm (z) \in \mathbb{Y}$ of finite type is reduced to exactly one pair in Table {\rm\ref{table:finie type}} by
 permuting the indices, changing sign, and changing slices $($see Section {\rm\ref{section:change of slice})}, if necessary.
 \end{enumerate}
\end{Theorem}

The claim (1) can be proved by concrete calculation in a suitable universal algebra
since~$A_\pm (z)$ in Table \ref{table:finie type} is concrete.
We, however, give another proof involving cluster algebras.
We give a~quiver and a sequence of mutations for each $A_\pm (z)$ in Table \ref{table:finie type}
that yields the $Y$-system as the exchange relation of coefficients
in the cluster algebra. See Table \ref{table:quivers} for quivers and mutations.
We can verify that some iteration of this sequence of mutations, as well as its inverse, is a~reddening sequence (see Theorem \ref{theorem:reddening}).
Thanks to the deep results in the theory of cluster algebras,
this property is enough to imply the periodicity (see Proposition \ref{prop:fin type of reddening}).
The numbers $h_\pm$ in Table~\ref{table:finie type} are the length of reddening sequences
in positive and negative directions, respectively.
This verification of the periodicity is interesting not only because it is computationally more efficient,
but also because it leads to nontrivial dilogarithm identities associated with
Donaldson--Thomas invariants (see Corollary \ref{cor:q-dilog identity}).

The claim (2) is proved in Section \ref{section:proof of classification} by the following steps:
\begin{enumerate}[label=Step \arabic*., leftmargin=*, itemindent=3em]\itemsep=0pt
 \item We recall the result in \cite{mizuno2020difference} that
 asserts that $A_\pm(1)$ satisfies a certain positivity,
 which in particular implies that $\tr A_\pm (1)$ and $\det A_\pm(1)$
 are positive. This allows us to significantly reduce the candidates for finite type $A_\pm(z)$.
 \item For a fixed $A_+(1)$ in the candidates obtained in Step 1, we search for $A_-(1)$
 satisfying the symplectic property \eqref{eq:symplectic property} at $z=1$.
 \item During the search in Step 2, we discard the pair $A_\pm (1)$ that cannot be endowed with
 the parameter $z$ (see Lemmas~\ref{lemma:ban list} and~\ref{lemma:ban list indec}).
 \item At this point, we have six candidates up to a permutation of the indices
 and a change of sign. For each $A_\pm(1)$ in the six candidates, we try to endow with the
 parameter $z$. It turns out that this is possible for all the six candidates.
 We give all possible $A_\pm(z)$ in Lemmas~\ref{lemma:Apm(z) 1 2 3}--\ref{lemma:Apm(z) 6}.
 \item We finally check that each remaining candidate reduces to one of $A_\pm (z)$ in
 Table~\ref{table:finie type} by change of slices.
\end{enumerate}

\begin{Remark}
 In fact, we can see that six pairs in Table \ref{table:finie type} are not related to each other
 by permuting the indices, changing sign, and changing slices.
 It can be verified by comparing the mutation sequences of quivers associated with these pairs,
 which are given in Table \ref{table:quivers}.
 See Remark \ref{remark:quiver is different}.
\end{Remark}

\begin{Remark}
 Most of the $Y$-systems obtained from Table \ref{table:finie type} are
 already known in the literature. $\mathrm{(1)}^{\mathrm{op}}$ and $\mathrm{(6)}^{\mathrm{op}}$
 are Zamolodchikov's $Y$-system of type $A_2$ \cite{Zamo} and $T_2$ (``tadpole'') \cite{RVT}, respectively.
 $\mathrm{(2)}^{\mathrm{op}}$ is the reduced sine-Gordon $Y$-system associated with the continued fraction~${3/4 = [1,3] = 1/(1+1/3)}$,
 and $\mathrm{(5)}^{\mathrm{op}}$ with $z$ replaced by $z^2$ is the reduced sine-Gordon $Y$-system associated with $3/5 = [1,1,2] = 1/(1+1/(1+1/2))$ \cite{Tateo} .
 (4) is the ``half'' of the $Y$-system associated with the pair $(A_2, A_2)$ \cite{RVT}.
 (3) appears to be new
 \begin{gather*}
 Y_1 (u) Y_1 (u - 2) = \frac{1}{1 + Y_2(u-1)^{-1}}, \\
 Y_2 (u) Y_2 (u - 10) =
 \frac{(1 + Y_1(u-3))(1 + Y_1(u-7))}
 {\bigl(1 + Y_1(u-1)^{-1}\bigr)\bigl(1 + Y_1(u-5)^{-1}\bigr)\bigl(1 + Y_1(u-9)^{-1}\bigr)}
 \end{gather*}
 although it is implicitly given in the author's previous work \cite[Table 2]{mizuno2020difference}.
\end{Remark}

\begin{Remark}
 The pair $A_\pm (z) \in \mathbb{Y}$ is called the \emph{T-datum} in \cite{mizuno2020difference} since
 it describes the T-systems, which is a companion to the $Y$-systems. We do not use this term since
 we only consider the $Y$-systems in this paper.
 Moreover, the definition of the T-datum in \cite{mizuno2020difference} allows to have a non-diagonal $N_0$ and
 have a nontrivial symmetrizer $D$, which is more general than the definition in this paper.
 See also Section \ref{section:comments} for the $Y$-systems involving nontrivial symmetrizers.
\end{Remark}

\begin{Remark}
 There is another expression of the $Y$-system using a pair of matrices $A_\pm (z)$ directly.
 Let $A_\pm (z) \in \mathbb{Y}_0$, and define $a_{ij;p} \in \mathbb{Z}$ by
 \begin{equation*}
 A_\pm (z) = \biggl( \sum_{p \in \mathbb{N}} a_{ij;p}^\pm z^p \biggr)_{i,j \in I}.
 \end{equation*}
 Let \smash{$\bigl(P_i^\pm (u)\bigr)_{i \in I, u \in \mathbb{Z}}$} be a family of elements in
 a multiplicative abelian group $\mathbb{P}$.
 We say that $\bigl(P_i^\pm (u)\bigr)$ \emph{satisfies the multiplicative $Y$-system} associated with $A_\pm (z)$ if
 \begin{align}\label{eq:y-system mul}
 \prod_{j \in I} \prod_{p \in \mathbb{N}} P_j^+ (u-p)^{a_{ij;p}^+}
 = \prod_{j \in I} \prod_{p \in \mathbb{N}} P_j^- (u-p)^{a_{ij;p}^-}
 \end{align}
 for any $i$, $u$ (schematically, ``$A_+(z) \cdot \log P^+ = A_-(z) \cdot \log P^-$''
 under the action $z\colon u \mapsto u-1$).
 The solution $\bigl(P_i^\pm (u)\bigr)$ is called \emph{normalized} if
 $\mathbb{P}$ is endowed with a semifield structure, and
$
 P_i^+(u) + P_i^-(u) = 1
$
 for any $i$, $u$.
 We have a one-to-one correspondence between solutions of the $Y$-system~\eqref{eq:y-system} and
 normalized solutions of the multiplicative $Y$-system \eqref{eq:y-system mul}.
 The correspondence is given by
 \begin{equation*}
 Y_i (u) \mapsto \frac{P_i^+(u)}{P_i^-(u)},\qquad
 P_i^+ (u) \mapsto \frac{Y_i(u)}{1+Y_i(u)},\qquad
 P_i^- (u) \mapsto \frac{1}{1+Y_i(u)}.
 \end{equation*}
 In the setting of cluster algebras, this correspondence is nothing but
 the normalization of the coefficients described by Fomin and Zelevinsky \cite[Section 5]{FZ1}.
\end{Remark}

\subsection{Relation to Nahm sums}\label{section:nahm sum}
Consider the $q$-series defined by
\begin{equation}\label{eq:RR sum}
 G(q) = \sum_{n=0}^\infty \frac{q^{n^2}}{(q)_n} ,\qquad
 H(q) = \sum_{n=0}^\infty \frac{q^{n^2+n}}{(q)_n},
\end{equation}
where $(q)_n \coloneqq (1-q) \bigl(1-q^2\bigr) \cdots (1-q^n)$ is the $q$-Pochhammer symbol.
The famous Rogers--Ramanujan identities express these $q$-series as the following infinite products:
\begin{equation*}
 G(q) = \prod_{n \equiv \pm 1 \bmod 5} \frac{1}{1 - q^n},\qquad
 H(q) = \prod_{n \equiv \pm 2 \bmod 5} \frac{1}{1 - q^n} .
\end{equation*}
These expressions in particular implies that $q^{-1/60} G(q)$ and $q^{11/60} H(q)$ are
modular functions on some finite index subgroup of $\SL (2, \mathbb{Z})$.
In fact, it is a rare case that an infinite sum of the form \eqref{eq:RR sum} is modular. It is known that the $q$-series
\begin{equation}\label{eq:Nahm sum 1}
 \sum_{n = 0}^\infty \frac{q^{ \frac{1}{2} a n^2 + b n + c}}{(q)_{n}}
\end{equation}
with $a,b,c \in \mathbb{Q}$ is modular only if $a = 1/2$, $1$, or $2$ \cite{Zagier}.

Nahm \cite{Nahm} considered higher rank generalization of \eqref{eq:Nahm sum 1},
which we call the \emph{Nahm sum}.
Let~$I$ be a finite set, and suppose that $A \in \mathbb{Q}^{I \times I}$ is a symmetric positive definite matrix,
$B \in \mathbb{Q}^I$ is a~vector,
and $C \in \mathbb{Q}$ is a scalar.
The Nahm sum is the $q$-series defined by
\begin{equation*}
 f_{A, B, C} (q) \coloneqq \sum_{n \in \mathbb{N}^I}
 \frac{q^{ \frac{1}{2} n^\mathsf{T} A n + n^\mathsf{T} B +C}}{ \prod_a (q)_{n_i}}.
\end{equation*}
When $\lvert I \rvert \geq 2$, it is not well understood when $f_{A,B,C}(q)$ is modular.
Nahm gave a conjecture providing a criterion on the modularity of $f_{A, B, C} (q)$
in terms of torsion elements in the Bloch group \cite{Nahm, Zagier}.
See \cite{calegari2017bloch, VlasenkoZwegers} for the development of this conjecture.

Nahm used Zamolodchikov's periodicity to provide an evidence of the conjecture.
In fact, there is a natural way to give a candidate of modular Nahm sums from finite
type $A_\pm (z) \in \mathbb{Y}$ in general.
Precisely, the matrix $K \coloneqq A_{+}(1)^{-1} A_{-}(1)$ is
always symmetric and positive definite for finite type $A_\pm (z) \in \mathbb{Y}$,
and it is conjectured that it gives a modular Nahm sum $f_{K, 0, C} (q)$ for some $C$ \cite{mizuno2020difference}.
(This construction is essentially the same as that in \cite{KatoTerashima}, except that they did not prove
that $K$ is symmetric and positive definite. A special case can also be found in \cite{Lee}.)
We~note that the symplectic property \eqref{eq:symplectic property} at $z=1$ plays an important role
here since it implies that $K$ is symmetric.
On the other hand, the positive definiteness is related to the periodicity of the $Y$-system.

\begin{table}[t]
 \centering
 \begin{tabular}{Sl | Sl | Sc | Sl | Sl | Sl | Sl | Sc | Sl} \hline
 $A_\pm (z)$ & $K$ & $-24C$ & RR & &
 $A_\pm (z)$ & $K$ & $-24C$ & RR \\ \hline
 (1) &
 $\begin{pmatrix}
 4/3 & 2/3 \\
 2/3 & 4/3
 \end{pmatrix}$ &
 $\dfrac{4}{5}$ &
 \cite{CherednikFeigin} & &
 $\mathrm{(1)}^{\mathrm{op}}$ &
 $\begin{pmatrix}
 1 & - 1/2 \\
 - 1/2 & 1
 \end{pmatrix}$ &
 $\dfrac{6}{5}$ &
 \cite{VlasenkoZwegers} \\
 (2) &
 $\begin{pmatrix}
 3/2 & 1 \\
 1 & 2
 \end{pmatrix}$ &
 $\dfrac{5}{7}$ &
 \cite{wang2022rank2} & &
 $\mathrm{(2)}^{\mathrm{op}}$ &
 $\begin{pmatrix}
 1 & -1/2 \\
 -1/2 & 3/4
 \end{pmatrix}$ &
 $\dfrac{9}{7}$ &
 \cite{wang2022rank2} \\
 (3), (6) &
 $\begin{pmatrix}
 2 & 2 \\
 2 & 4
 \end{pmatrix}$ &
 $\dfrac{4}{7}$ &
 \cite{Andrews74} & &
 $\mathrm{(3)}^{\mathrm{op}}$,
 $\mathrm{(6)}^{\mathrm{op}}$ &
 $\begin{pmatrix}
 1 & -1/2 \\
 -1/2 & 1/2
 \end{pmatrix}$ &
 $\dfrac{10}{7}$ &
 \cite{wang2022rank2} \\
 (4) &
 $\begin{pmatrix}
 2/3 & 1/3 \\
 1/3 & 2/3
 \end{pmatrix}$ &
 $1$ &
 \cite{Zagier} & &
 $\mathrm{(4)}^{\mathrm{op}}$ &
 $\begin{pmatrix}
 2 & -1 \\
 -1 & 2
 \end{pmatrix}$ &
 $1$ &
 \cite{Zagier} \\
 (5) &
 $\begin{pmatrix}
 1 & 1 \\
 1 & 2
 \end{pmatrix}$ &
 $\dfrac{3}{4}$ &
 \cite{wang2022rank2} & &
 $\mathrm{(5)}^{\mathrm{op}}$ &
 $\begin{pmatrix}
 2 & -1 \\
 -1 & 1
 \end{pmatrix}$ &
 $\dfrac{5}{4}$ &
 \cite{CalinescuMilasPenn16} \\
 \hline
 \end{tabular}
 \caption{The list of the matrix $K = A_+(1)^{-1} A_- (1)$. The Nahm sum $f_{K,0,C}(q)$ is modular,
 which can be proved by using Rogers--Ramanujan type identities (RR for short) given in the references
 in the table.} \label{table:Nahm sum}
\end{table}

Based on our classification, we find the following statement.
\begin{Theorem}
 Suppose that $I= \{1, 2\}$.
 The Nahm sum $f_{K, 0, C} (q)$ is modular for any finite type $A_\pm (z) \in \mathbb{Y}$,
 where $C$ is given in Table~{\rm\ref{table:Nahm sum}}.
\end{Theorem}
In fact, every $K$ from finite type $A_\pm(z)$ is included in the Zagier's list \cite[Table 2]{Zagier}
for rank~$2$ candidates of modular Nahm sums.
There are Rogers--Ramanujan type identities that enable us to write each Nahm sum in the list in terms of theta functions.
The proof of the desired identities was partially given in \cite{Andrews74,CalinescuMilasPenn16,CherednikFeigin, VlasenkoZwegers,Zagier},
and was recently completed by Wang \cite{wang2022rank2} except for one candidate that does not appears in our construction from $Y$-systems.
See Table \ref{table:Nahm sum}.

\begin{Example}
 For a matrix $K =
 \bigl(\begin{smallmatrix}
 2 & 2 \\
 2 & 4
 \end{smallmatrix}\bigr)$,
 we have
 \begin{equation*}
 f_{K, (0, 0), -1/42} (q) = \frac{\theta_{7, 1}(\tau)}{\eta(\tau)}, \qquad
 f_{K, (0, 1), 5/42} (q) = \frac{\theta_{7, 2}(\tau)}{\eta(\tau)}, \qquad
 f_{K, (1, 2), 17/42} (q) = \frac{\theta_{7, 3}(\tau)}{\eta(\tau)},
 \end{equation*}
 by the Andrews--Gordon identities \cite{Andrews74}
 and the Jacobi's triple product identity,
 where
 \begin{equation*}
 \eta (\tau) \coloneqq q^{1/24} \prod_{n=1}^\infty (1 - q^n), \qquad
 \theta_{7,j} (\tau) \coloneqq \sum_{n \equiv 2 j -1 \bmod{14}} (-1)^{[n/14]}
 q^{n^2 / 56}, \qquad
 \end{equation*}
 with $q = {\rm e}^{2 \pi {\rm i} \tau}$.
 The Dedekind eta function $\eta(\tau)$ and the theta functions $\theta_{7,j} (\tau)$
 are modular forms of weight $1/2$,
 and their modular transformation formulas imply that the vector-valued function
 $f = (f_{K,0, -1/42} f_{K, (0,1), 5/42} f_{K, (1,2), 17/42})^{\mathsf{T}}$
 is a vector-valued modular function satisfying
 \begin{equation*}
 f(-1/\tau) =
 \frac{2}{\sqrt{7}}
 \begin{pmatrix}
 \sin \frac{3 \pi}{7} & \sin \frac{2 \pi}{7} & \sin \frac{\pi}{7} \\
 \sin \frac{2 \pi}{7} & -\sin \frac{\pi}{7} & -\sin \frac{3 \pi}{7} \\
 \sin \frac{\pi}{7} & -\sin \frac{3 \pi}{7} & \sin \frac{2 \pi}{7}
 \end{pmatrix}
 f(\tau).
 \end{equation*}
\end{Example}

\begin{Remark}
 We can define the refinement \smash{$f_{A_\pm (z)}^{(s)} (q)$} of the Nahm sum $f_{K,0,0} (q)$,
 which is parametrized by $s \in H$ for an abelian group $H$ of order $\det A_+ (1)$
 such that~it reduces~to~the~original one by taking summation \cite[Definition 5.12]{mizuno2020difference}
 \begin{equation*}
 f_{K,0,0} (q) = \sum_{s \in H} f_{A_\pm (z)}^{(s)} (q).
 \end{equation*}
 It is conjectured that each \smash{$f_{A_\pm (z)}^{(s)} (q)$} is already
 modular after multiplying $q^C$ for some $C$.
 We note that the symplectic property \eqref{eq:symplectic property}
 at $z=1$ again plays an important role in the definition of the refinement.
 We will discuss this refinement for rank 2 case in more detail elsewhere.
 We remark that similar refinement also appears in the context of
 $3$-dimensional quantum topology~\cite[Section 6.3]{garoufalidis2021knots}.
\end{Remark}

\subsection{Remarks on higher rank and skew-symmetrizable case}
\label{section:comments}
We have seen that the following properties hold for rank $2$ case:
\begin{enumerate}[label=(P\arabic*)]\itemsep=0pt
 \item \label{item:P1} We have reddening sequences in both positive and negative directions.
 \item \label{item:P2} The map $A_\pm(z) \mapsto A_+ (1)^{-1} A_- (1)$ gives modular Nahm sums.
\end{enumerate}
We expect that the properties \ref{item:P1} and \ref{item:P2} also hold for any finite
type $A_\pm (z) \in \mathbb{Y}$ of general rank. The following are some known examples:
\begin{itemize}\itemsep=0pt
 \item For the $Y$-system associated with the untwisted quantum affine algebras $U_q\bigl(X_r^{(1)}\bigr)$
 with level $\ell$ restriction \cite{KunibaNakanishi92},
 \ref{item:P1} holds with $h_+ = t \ell$ and $h_- = t \cdot (\text{dual Coxeter number of $X_r$})$
 where $t = 1$, $2$, or $3$ is the multiplicity in the Dynkin diagram of $X_r$ \cite{IIKKNa,IIKKNb}, and
 \ref{item:P2} holds under the assumption \cite[Conjecture 5.3]{HKOTT} by the
 result of Kac and Peterson \cite{KacPet}.
 \item For the $Y$-system associated with a pair of finite type simply laced Dynkin type $(X_r, X'_{r'})$ \cite{RVT},
 \ref{item:P1} holds with $h_+ = (\text{Coxeter number of $X_r$})$ and $h_- = (\text{Coxeter number of $X'_{r'}$})$ \cite{Keller2011, Keller}.
 \item For the (reduced) sine-Gordon $Y$-system associated with the continued fraction
 $p/q = [n_F, \dots, n_1] = 1/ (n_F + 1/ (\cdots + 1/ n_1))$ \cite{NakanishiStella, Tateo},
 \ref{item:P1} appears to hold with $h_+ = 2 p$ and $h_- = 2 q$.
 \item For the $Y$-system associated with an admissible ADE bigraph $(\Gamma, \Delta)$ \cite{GalashinPylyavskyy},
 \ref{item:P1} appears to hold with $h_+ = (\text{Coxeter number of $\Gamma$})$ and
 $h_- = (\text{Coxeter number of $\Delta$})$.
\end{itemize}

Moreover, we can consider $Y$-systems associated with skew-symmetrizable cluster algebras rather than skew-symmetric ones
discussed in this paper. In this case, the symplectic property~\eqref{eq:symplectic property} becomes
\begin{equation*}
 A_+ (z) D A_-\bigl(z^{-1}\bigr)^{\mathsf{T}} = A_- (z) D A_+ \bigl(z^{-1}\bigr)^{\mathsf{T}},
\end{equation*}
where $D$ is a diagonal matrix called symmetrizer \cite{mizuno2020difference}.
We also expect that the properties \ref{item:P1} and \ref{item:P2} also hold for skew-symmetrizable case.
See \cite[Definition 5.12]{mizuno2020difference} for the definition of the Nahm sum in skew-symmetrizable case.

\section[Y-systems and cluster algebras]{$\boldsymbol{ Y}$-systems and cluster algebras}
\subsection{Preliminaries on cluster algebras}
\label{section:Y from cluster}
In this paper, a \emph{semifield} is a multiplicative abelian group equipped with an addition that is commutative,
associative, and distributive with respect to the multiplication.

\begin{Definition}
 Let $I$ be a set.
 The set of all nonzero rational functions in the variables $\boldsymbol{y} = (y_i)_{i \in I}$ with
 natural number coefficients is a semifield with respect to the usual addition and multiplication.
 This semifield is called the \emph{universal semifield}, and denoted by $\mathbb{Q}_{>0} (\boldsymbol{y})$.
 We have a canonical bijection
 $\Hom_{\mathrm{semifield}} (\mathbb{Q}_{>0} (\boldsymbol{y}), \mathbb{P}) \cong \Hom_{\mathrm{set}} (I, \mathbb{P})$
 for any set $I$ and semifield $\mathbb{P}$.
\end{Definition}

\begin{Definition}
 Let $I$ be a set.
 The \emph{tropical semifield} $\Trop(\boldsymbol{y})$ is the multiplicative free abelian group
 generated by the variables $\boldsymbol{y} = (y_i)_{i \in I}$ equipped with the addition defined by
 \begin{equation*}
 \prod_i y_i^{a_i} +
 \prod_i y_i^{b_i} =
 \prod_i y_i^{\min(a_i,b_i)}.
 \end{equation*}
\end{Definition}

Let $I$ be a finite set and $\mathbb{P}$ be a semifield.
A \emph{Y-seed} is a pair $(B, \boldsymbol{y})$ where $B = (B_{ij})_{i,j \in I}$ is a skew-symmetric integer matrix and
$\boldsymbol{y} = (y_i)_{i \in I}$ is a tuple of elements in $\mathbb{P}$.
We sometimes represent $B$ as the quiver whose signed adjacency matrix is $B$.
For a Y-seed $(B, \boldsymbol{y})$ and $k \in I$, the \emph{mutation} in direction $k$ transforms $(B, \boldsymbol{y})$ into the
new Y-seed $\mu_k(B, \boldsymbol{y}) = (B', \boldsymbol{y}')$ given by
\begin{gather*}
 B'_{ij} \coloneqq
 \begin{cases}
 -B_{ij} & \text{if $i=k$ or $j=k$}, \\
 B_{ij} +
 \maxzero{- B_{ik}} B_{kj} +
 B_{ik} \maxzero{ B_{kj}} & \text{otherwise},
 \end{cases} \\
 y'_i \coloneqq
 \begin{cases}
 y_k & \text{if $i=k$}, \\
 y_i y_k^{\maxzero{B_{ki}}} (1 + y_k)^{-B_{ki}}
 & \text{otherwise}.
 \end{cases}
\end{gather*}
A mutation is involutive, that is, $\mu_k (B, \boldsymbol{y}) = (B', \boldsymbol{y}')$ implies
$(B, \boldsymbol{y}) = \mu_k (B', \boldsymbol{y}')$.
We have the commutativity
\begin{equation}\label{eq:mutation comm}
 \mu_i \mu_j = \mu_j \mu_i \qquad \text{if}\quad B_{ij} = 0,
\end{equation}
which allows us to write $\mu_{\boldsymbol{i}}$ for a set $\boldsymbol{i} \subseteq I$
such that $B_{ij} = 0$ for any $i,j \in \boldsymbol{i}$
to mean the successive mutations along arbitrarily chosen order on $\boldsymbol{i}$.

For a Y-seed $(B, \boldsymbol{y})$ and a bijection $\nu \colon I \to I$, we define a new Y-seed
$\nu(B, \boldsymbol{y}) = (B', \boldsymbol{y}')$ by~${B'_{\nu(i)\nu(j)} \coloneqq B_{ij}}$ and $y'_{\nu(i)} \coloneqq y_i$.

\subsection[Solving Y-systems by cluster algebras]{Solving $\boldsymbol{ Y}$-systems by cluster algebras}
Let $A_\pm (z) \in \mathbb{Y}$.
We will construct a solution
of the $Y$-system associated with $A_\pm (z)$ based on~\cite[Section 3.3]{mizuno2020difference}.
We first define a subset $R \subseteq I \times \mathbb{Z}$ by
$
 R\coloneqq \{(i, u) \in I \times \mathbb{Z} \mid 0 \leq u < r_i \}$,
and define a skew-symmetric $R \times R$ integer matrix $B$ by
\begin{equation}\label{eq:def of B from cap}
 B_{(i, p)(j, q)} = - n_{ij;p-q} + n_{ji;q-p} +\sum_{k \in I} \sum_{v=0}^{\min(p,q)}
 \bigl( n_{ik;p-v}^{+} {n}_{jk;q-v}^{-}
 - n_{ik;p-v}^{-} {n}_{jk;q-v}^{+} \bigr),
\end{equation}
where we understand $n_{ij;p} = 0$ if $p < 0$.
We then define $\boldsymbol{i} \coloneqq \{ (i, u) \mid u=0 \} \subseteq R$.
We also define a bijection $\nu \colon R \to R$ by
\[
	\nu (i, p) =
	\begin{cases}
		(i, p - 1) &\text{if $p > 0$,}\\
		(i, r_i) &\text{if $p=0$.}
	\end{cases}
\]
Then the symplectic property \eqref{eq:symplectic property} ensures that $\nu(\mu_{\boldsymbol{i}} (B)) = B$ \cite[Lemma 3.16]{mizuno2020difference}.
We finally define a sequence of Y-seeds
\begin{align}\label{eq:seq of y universal}
 \cdots \to (B, \boldsymbol{y}(-1)) \to (B, \boldsymbol{y}(0)) \to (B, \boldsymbol{y}(1)) \to \cdots
\end{align}
in $\mathbb{Q}_{> 0}(\boldsymbol{y})$ by
$\boldsymbol{y}(0) \coloneqq \boldsymbol{y}$ and
$(B, \boldsymbol{y}(u+1)) = \nu (\mu_{\boldsymbol{i}} (B, \boldsymbol{y}(u)))$.
The sequence \eqref{eq:seq of y universal} gives a solution of the $Y$-system.

\begin{Lemma}[{\cite[Theorem 3.13]{mizuno2020difference}}]\label{lemma:universal solution}
 $(y_{i, 0} (u))_{i \in I, u \in \mathbb{Z}}$ satisfies the $Y$-system associated with $A_\pm(z)$.
\end{Lemma}

This solution is universal in the following sense.
\begin{Lemma}[{\cite[Theorem 3.19]{mizuno2020difference}}]
 Suppose that a family $(Y_i(u))_{i \in I, u \in \mathbb{Z}}$ satisfies the $Y$-system associated with $A_\pm(z)$.
 Define a semifield homomorphism $f \colon \mathbb{Q}_{>0} (\boldsymbol{y}) \to \mathbb{P}$ by
 \[
		f(y_{i, p}) \coloneqq
 Y_i(p) \prod_{j \in I} \prod_{q=0}^{p} Y_j (p-q)^{- \maxzero{ n_{ij;q} }} (1 + Y_j (p-q) )^{n_{ij;q}}.
	\]
 Then $f(y_{i, 0} (u)) = Y_i(u)$ for any $i$, $u$.
\end{Lemma}

\begin{Corollary}\label{cor:finite iff periodic for universal y}
 $A_{\pm} (z) \in \mathbb{Y}$ is of finite type if and only if there are different integers $u$, $v$ such that
 $\boldsymbol{y}(u) = \boldsymbol{y} (v)$ in \eqref{eq:seq of y universal}.
\end{Corollary}

\subsection{Periodicity and reddening sequences}
Similarly to \eqref{eq:seq of y universal}, we define a sequence of Y-seeds
\begin{align}\label{eq:seq of y tropical}
 \cdots \to (B, \boldsymbol{y}(-1)) \to (B, \boldsymbol{y}(0)) \to (B, \boldsymbol{y}(1)) \to \cdots
\end{align}
by the same formulas but now in $\Trop(\boldsymbol{y})$ rather than $\mathbb{Q}(\boldsymbol{y})$.

\begin{Definition}
 We say that the $Y$-system associated with $A_\pm (z) \in \mathbb{Y}$ is \emph{positive $($resp.\ negative$)$ reddening} if
 there is a positive integer $u$ such that
 all the exponents in $y_i(u)$ (resp.\ $y_i(-u)$) in~\eqref{eq:seq of y tropical} are nonpositive for any $i$.
 We denote by $h_+$ (resp.\ $h_-$) the least such positive integer $u$.
\end{Definition}

Equivalently, the $Y$-system is positive (resp.\ negative) reddening if and only if all the entries in the
C-matrix associated with the sequence of mutations
$(B, \boldsymbol{y}(0)) \to (B, \boldsymbol{y}(u))$
(resp.\ $(B, \boldsymbol{y}(0)) \to (B, \boldsymbol{y}(-u))$)
are nonpositive for some $u>0$.

\begin{Proposition}\label{prop:fin type of reddening}
 Suppose that the $Y$-system associated with $A_\pm (z)$ is positive and negative reddening.
 Then $A_\pm (z)$ is of finite type.
\end{Proposition}
\begin{proof}
 We verify the equivalent condition in Corollary \ref{cor:finite iff periodic for universal y}.
 By \cite[Proposition 2.10]{BDP},
 there are bijections $\sigma, \sigma' \colon R \to R$ such that
 \smash{$y_i(h_+) = y_{\sigma(i)}^{-1}$} and \smash{$y_i(-h_-) = y_{\sigma'(i)}^{-1}$} for any $i$
 (in other words, the C-matrices associated with them are the minus of permutation matrices).
 Now the claim follows from the separation formula for $y$-variables \cite[Proposition 3.13]{FZ4} and
 the result on C-matrices shown by Cao, Huang, and Li \cite[Theorem 2.5]{CHL20}.
 See also \cite[Theorem 5.2]{nakanishi2019synchronicity} for the corresponding statement dealing with permutations
 that is actually suitable here.
\end{proof}

\begin{table}[t]
 \begin{center}
 \begin{tabular}{ Sc | Sl } \hline
 Quiver $B$ & $A_\pm(z)$ \\ \hline
 \begin{tikzpicture}[scale=1.6, baseline=(00.base)]
 \node (00) at (0,0) {$(1,0)$};
 \draw (00) ++ (1,0) node (11) {$(2,1)$};
 \draw [quiver] (00) -- (11);
 \newcommand{\x}{2.5}
 \node (01) at (\x,0) {$(1,1)$};
 \draw (01) ++ (1,0) node (10) {$(2,0)$};
 \draw [quiver] (10) -- (01);
 \end{tikzpicture} & (1) \\ \hline
 \begin{tikzpicture}[scale=1.5, baseline=(00.base)]
 \node (00) at (0,0) {$(1,0)$};
 \draw (00) ++ (0:1) node (11) {$(2,1)$};
 \draw (00) ++ (120:1) node (13) {$(2,3)$};
 \draw (00) ++ (240:1) node (15) {$(2,5)$};
 \draw [quiver] (00) -- (11);
 \draw [quiver] (13) -- (00);
 \draw [quiver] (00) -- (15);
 \draw [quiver] (11) -- (13);
 \newcommand{\x}{2.5}
 \node (01) at (\x,0) {$(1,1)$};
 \draw (01) ++ (0:1) node (10) {$(2,0)$};
 \draw (01) ++ (120:1) node (12) {$(2,2)$};
 \draw (01) ++ (240:1) node (14) {$(2,4)$};
 \draw [quiver] (10) -- (01);
 \draw [quiver] (01) -- (12);
 \draw [quiver] (14) -- (01);
 \draw [quiver] (12) -- (14);
 \end{tikzpicture} & (2) \\ \hline
 \begin{tikzpicture}[scale=1.5, baseline=(00.base)]
 \node (00) at (0,0) {$(1,0)$};
 \draw (00) ++ (0:1) node (11) {$(2,1)$};
 \draw (00) ++ (360/5:1) node (13) {$(2,3)$};
 \draw (00) ++ (360*2/5:1) node (15) {$(2,5)$};
 \draw (00) ++ (360*3/5:1) node (17) {$(2,7)$};
 \draw (00) ++ (360*4/5:1) node (19) {$(2,9)$};
 \draw [quiver] (00) -- (11);
 \draw [quiver] (13) -- (00);
 \draw [quiver] (00) -- (15);
 \draw [quiver] (17) -- (00);
 \draw [quiver] (00) -- (19);
 \draw [quiver] (11) -- (13);
 \draw [quiver] (11) -- (17);
 \draw [quiver] (15) -- (17);
 \newcommand{\x}{2.5}
 \node (01) at (\x,0) {$(1,1)$};
 \draw (01) ++ (0:1) node (10) {$(2,0)$};
 \draw (01) ++ (360/5:1) node (12) {$(2,2)$};
 \draw (01) ++ (360*2/5:1) node (14) {$(2,4)$};
 \draw (01) ++ (360*3/5:1) node (16) {$(2,6)$};
 \draw (01) ++ (360*4/5:1) node (18) {$(2,8)$};
 \draw [quiver] (10) -- (01);
 \draw [quiver] (01) -- (12);
 \draw [quiver] (14) -- (01);
 \draw [quiver] (01) -- (16);
 \draw [quiver] (18) -- (01);
 \draw [quiver] (12) -- (14);
 \draw [quiver] (12) -- (18);
 \draw [quiver] (16) -- (18);
 \end{tikzpicture} & (3) \\ \hline
 \begin{tikzpicture}[scale=1.5, baseline={(0,1/2)}]
 \node (00) at (0,0) {$(1,0)$};
 \draw (00) ++ (1,0) node (11) {$(2,1)$};
 \draw (11) ++ (0,1) node (10) {$(2,0)$};
 \draw (10) ++ (-1,0) node (01) {$(1,1)$};
 \draw [quiver] (00) -- (11);
 \draw [quiver] (10) -- (11);
 \draw [quiver] (10) -- (01);
 \draw [quiver] (01) -- (00);
 \end{tikzpicture} & (4) \\ \hline
 \begin{tikzpicture}[scale=1.6, baseline={(0,0.7/2)}]
 \node (10) at (0,0) {$(2,0)$};
 \draw (10) ++ (1,0) node (01) {$(1,1)$};
 \draw (01) ++ (1,0) node (00) {$(1,0)$};
 \draw (00) ++ (1,0) node (12) {$(2,2)$};
 \draw (00) ++ (-1/2,0.7) node (11) {$(2,1)$};
 \draw [quiver] (10) -- (01);
 \draw [quiver] (01) -- (00);
 \draw [quiver] (00) -- (12);
 \draw [quiver] (00) -- (11);
 \draw [quiver] (11) -- (01);
 \end{tikzpicture} & (5) \\ \hline
 \begin{tikzpicture}[scale=1.6, baseline=(00.base)]
 \node (00) at (0,0) {$(1,0)$};
 \draw (00) ++ (1,0) node (11) {$(2,1)$};
 \draw (11) ++ (1,0) node (10) {$(2,0)$};
 \draw (10) ++ (1,0) node (01) {$(1,1)$};
 \draw [quiver] (00) -- (11);
 \draw [quiver] (10) -- (11);
 \draw [quiver] (10) -- (01);
 \end{tikzpicture} & (6) \\
 \hline
 \end{tabular}
 \end{center}
 \caption{Quivers associated with $A_\pm(z)$ in Table \ref{table:finie type}.
 Each quiver is preserved by the mutation at $({*}, 0)$ followed by the permutation $(i, p) \mapsto (i, p-1)$ (the second argument
 is considered modulo $r_i$), which yields $Y$-system. For (1)--(3), this operation interchanges the connected components
 (see Section \ref{section:change of slice}).}
 \label{table:quivers}
\end{table}

\begin{Theorem}\label{theorem:reddening}
 The $Y$-system associated with each $A_\pm(z)$ in Table {\rm\ref{table:finie type}} is positive and negative reddening.
\end{Theorem}
\begin{proof}
 The quiver $B$ associated with $A_\pm(z)$ is given in Table \ref{table:quivers}.
 We can verify the assertion by concrete calculation on the quiver.
 The numbers $h_\pm$ are given in Table \ref{table:finie type}.
\end{proof}

Theorem \ref{theorem:main}\,(1) now follows from Proposition \ref{prop:fin type of reddening} and Theorem \ref{theorem:reddening}.

\begin{Remark}
 A connected component of each quiver in Table \ref{table:quivers} has the following cluster type:
\[
 (1) \ A_2, \quad (2) \ A_4, \quad (3) \ E_6, \quad (4) \ D_4, \quad (5) \ A_5, \quad (6) \ A_4.
\]
 These are of finite type in the sense of~\cite{FZ2}, which also implies Theorem \ref{theorem:main}\,(1).
 We remark, however, that this observation is somewhat misleading since the quiver associated with
 a periodic $Y$-system of general rank is typically of infinite type.
 It might be better to think that the appearance of only finite type quivers happens ``by chance''
 due to the smallness of~$2$, the rank of $Y$-systems considered in this paper.
\end{Remark}

Theorem \ref{theorem:reddening} also gives quantum dilogarithm identities
associated with Donaldson--Thomas invariants.
For any reddening sequence $\boldsymbol{i}$ starting from a quiver $B$, we can define
a quantity~$\mathbb{E}(\boldsymbol{i})$ by using the quantum dilogarithm.
We refer to \cite[Remark 6.6]{keller12} as the definition.
This quantity coincides with Kontsevich--Soibelman's refined
Donaldson--Thomas invariant associated with~$B$~\cite{keller12, Nagao10}.
In particular, $\mathbb{E}(\boldsymbol{i})$ does not depend on $\boldsymbol{i}$,
which gives the quantum dilogarithm identities. In our case, we have the following.

\begin{Corollary}\label{cor:q-dilog identity}
 For each $A_\pm (z)$ in Table {\rm\ref{table:finie type}},
 we have
 \begin{equation*}
 \mathbb{E}\bigl(\mu^{h_+}\bigr) = \mathbb{E}\bigl(\mu^{- h_-}\bigr),
 \end{equation*}
 where $\mu \coloneqq \nu \circ \mu_{\boldsymbol{i}}$ is the sequence of mutations $($together with the permutation$)$
 $(B, \boldsymbol{y}(0)) \to (B, \boldsymbol{y}(1))$ in \eqref{eq:seq of y universal}.
\end{Corollary}

For example, the pair (1) in Table \ref{table:finie type} yields
the famous pentagon identity of the quantum dilogarithm.

\section{Classification}
\subsection{Change of slices}\label{section:change of slice}
We need to introduce an appropriate equivalence relation on the set $\mathbb{Y}$,
which identifies essentially the same $Y$-systems.
Before we get into the definition, we will see a typical example.
Consider the following $Y$-system
\begin{gather}\label{eq:slice example 1}
 Y_1 (u) Y_1 (u - 2) = \bigl(1 + Y_2(u-1)^{-1}\bigr)^{-1},\qquad Y_2 (u) Y_2 (u - 2) = \bigl(1 + Y_1(u-1)^{-1}\bigr)^{-1},
\end{gather}
which corresponds to $A_\pm (z) \in \mathbb{Y}$ given by (1) in Table \ref{table:finie type}.
This system of equations are defined on the set $[1,2] \times \mathbb{Z}$, but
actually can be defined on each component of the following disjoint union
\begin{equation*}
 [1,2] \times \mathbb{Z} = \bigsqcup_{k=0}^1 \{ (i, u) \mid i-u \equiv k \bmod 2 \}.
\end{equation*}
We informally call the algebraic relation defined on each subset
the \emph{slice} of the whole $Y$-system.
If $(Y_i(u))$ is a solution of the $Y$-system for $i-u \equiv 0 \bmod 2$,
then $(Y_i(u+1))$ is a solution of the $Y$-system for $i-u \equiv 1 \bmod 2$.
Thus it is enough to consider only one slice when considering solutions.
Now we consider another $Y$-system
\begin{gather}\label{eq:slice example 2}
Y'_1 (u) Y'_1 (u - 3) = \bigl(1 + Y'_2(u-2)^{-1}\bigr)^{-1},\qquad
 Y'_2 (u) Y'_2 (u - 3) = \bigl(1 + Y'_1(u-1)^{-1}\bigr)^{-1}.
\end{gather}
which corresponds to $A'_\pm (z) \in \mathbb{Y}$ given by
\begin{align*}
 A'_+(z)\coloneqq
 \begin{pmatrix}
 1 + z^3 & -z^2 \\
 -z & 1 + z^3
 \end{pmatrix},\qquad
 A'_-(z) \coloneqq
 \begin{pmatrix}
 1 + z^3 & 0 \\
 0 & 1 + z^3
 \end{pmatrix}.
\end{align*}
The $Y$-system \eqref{eq:slice example 2} is decomposed into three slices
\begin{equation*}
 [1,2] \times \mathbb{Z} =
 \bigsqcup_{k=0}^2 \{ (i, u) \mid i - u \equiv k \bmod 3 \}.
\end{equation*}
We see that for any solution of \eqref{eq:slice example 1} for $i- u \equiv 0 \bmod 2$,
\begin{equation*}
 Y'_1(u) \coloneqq Y_1 \left(\frac{2}{3}u - \frac{1}{3}\right),\qquad
 Y'_2(u) \coloneqq Y_2 \left(\frac{2}{3}u\right)
\end{equation*}
is a solution of \eqref{eq:slice example 2} for $i - u \equiv 2 \bmod 3$.
We also obtain solutions for the other two slices by shifting $u$.
Conversely, any solution of \eqref{eq:slice example 1} is obtained from a solution of \eqref{eq:slice example 2}.
Therefore, it~is enough to consider one of the $Y$-systems \eqref{eq:slice example 1} and \eqref{eq:slice example 2}.
In particular, $A_\pm(z)$ is of finite type if and only if $A'_\pm(z)$ is.

Now we work in the general setting.
The idea is that each slice corresponds to each connected component of the quiver
associated with the matrix $B$ defined by \eqref{eq:def of B from cap}.
Let $A_\pm(z) \in \mathbb{Y}$, and assume that it is indecomposable.
By \cite[Proposition 3.24]{mizuno2020difference}, we have a decomposition of
the matrix $B$ and its index set $R$
\begin{align*}
 B = \bigoplus_{u=0}^{t-1} B(u), \qquad
 R = \bigsqcup_{u=0}^{t-1} R(u)
\end{align*}
such that each $B(u)$ is indecomposable and we have a cyclic sequence of mutations
\begin{align}\label{eq:cyclic connected component}
 B(0) \xrightarrow{\nu|_{R(0)} \circ \mu_{\boldsymbol{i}(0)}} B(1)
 \longrightarrow
 \cdots
 \longrightarrow
 B(t-1)
 \xrightarrow{\nu|_{R(t-1)} \circ \mu_{\boldsymbol{i}(t-1)}}
 B(0),
\end{align}
where $\boldsymbol{i}(u) \coloneqq \boldsymbol{i} \cap R(u)$.
We say that two pairs $A_\pm(z)$ and $A'_\pm (z)$ are \emph{related by change of slices}
if they yield the same cyclic sequence \eqref{eq:cyclic connected component} up to a change
of indices and the commutativity of mutations \eqref{eq:mutation comm}. (This commutativity is
already implicitly used to justify the notation $\mu_{\boldsymbol{i}(u)}$ as stated below \eqref{eq:mutation comm}.)

\begin{Example}
 The pairs $A_\pm(z)$ and $A'_\pm(z)$ associated with \eqref{eq:slice example 1}
 and \eqref{eq:slice example 2}, respectively, are related by change of slices.
 Indeed, we see that the sequence \eqref{eq:cyclic connected component} for \eqref{eq:slice example 1} is
 \begin{equation*}
 \begin{tikzpicture}[scale=1.5, baseline=(00.base)]
 \node (00) at (0,0) {$(1,0)$};
 \node (11) at (1,0) {$(2,1)$};
 \draw [quiver] (00) -- (11);
 \end{tikzpicture}
 \xrightarrow{\nu \circ \mu_{(0,0)}}
 \begin{tikzpicture}[scale=1.5, baseline=(01.base)]
 \node (01) at (0,0) {$(1,1)$};
 \node (10) at (1,0) {$(2,0)$};
 \draw [quiver] (10) -- (01);
 \end{tikzpicture}
 \xrightarrow{\nu \circ \mu_{(1,0)}}
 \begin{tikzpicture}[scale=1.5, baseline=(00.base)]
 \node (00) at (0,0) {$(1,0)$};
 \node (11) at (1,0) {$(2,1)$,};
 \draw [quiver] (00) -- (11);
 \end{tikzpicture}
 \end{equation*}
 whereas the sequence \eqref{eq:cyclic connected component} for \eqref{eq:slice example 2} is
 \begin{equation*}
 \begin{tikzpicture}[scale=1.5, baseline=(00.base)]
 \node (00) at (0,0) {$(1,0)$};
 \node (11) at (1,0) {$(2,1)$};
 \draw [quiver] (00) -- (11);
 \end{tikzpicture}
 \xrightarrow{\nu' \circ \mu_{(0,0)}}
 \begin{tikzpicture}[scale=1.5, baseline=(02.base)]
 \node (02) at (0,0) {$(1,2)$};
 \node (10) at (1,0) {$(2,0)$};
 \draw [quiver] (10) -- (02);
 \end{tikzpicture}
 \xrightarrow{\nu' \circ \mu_{(1,0)}}
 \begin{tikzpicture}[scale=1.5, baseline=(01.base)]
 \node (01) at (0,0) {$(1,1)$};
 \node (12) at (1,0) {$(2,2)$};
 \draw [quiver] (01) -- (12);
 \end{tikzpicture}
 \xrightarrow{\nu'}
 \begin{tikzpicture}[scale=1.5, baseline=(00.base)]
 \node (00) at (0,0) {$(1,0)$};
 \node (11) at (1,0) {$(2,1)$.};
 \draw [quiver] (00) -- (11);
 \end{tikzpicture}
 \end{equation*}
 These are the same sequence up to a change of indices.
\end{Example}

\begin{Remark}\label{remark:quiver is different}
 The mutation sequences associated with the quivers in Table \ref{table:quivers}
 are not related to each other by
 a change of indices and the commutativity of mutations \eqref{eq:mutation comm}.
\end{Remark}

\subsection{Proof of the classification}
\label{section:proof of classification}
In this section, we will prove Theorem \ref{theorem:main}\,(2). We first recall the following result.
\begin{Lemma}[{\cite[Theorem 5.5]{mizuno2020difference}}]
 \label{lemma:simultaneous positive}
 Let $A_\pm(z) \in \mathbb{Y}$. Assume that $A_\pm(z)$ is of finite type.
 Then there is a vector $v \in \mathbb{R}^I$ such that $v > 0$,
 $v A_+(1) > 0$, and $v A_-(1) > 0$.
 In particular, $\tr A_\pm (1) >0$ and $\det A_\pm (1) > 0$.
\end{Lemma}

Since the off-diagonal entries of $A_\pm(1)$ are nonpositive integers by
the definition of $A_\pm(1)$,
the condition $\det A_\pm(1) > 0$ in Lemma \ref{lemma:simultaneous positive} implies that
the product of diagonal entries of $A_\pm(1)$ is positive.
By combining the condition $\tr A_\pm(1) > 0$ in Lemma \ref{lemma:simultaneous positive},
we see that the diagonal entries of $A_\pm(1)$ are positive.
By the definition of $A_\pm(1)$, the diagonal entries of $A_\pm(1)$ are less than or equal to $2$,
and hence they are equal to $1$ or $2$.
Consequently, we see that $A_+(1)$ and~$A_-(1)$ are equal to one of the following matrices:
\begin{gather*}
 \begin{pmatrix}
 2 & -1 \\
 -1 & 2
 \end{pmatrix},\qquad
 \begin{pmatrix}
 2 & -1 \\
 -2 & 2
 \end{pmatrix},\qquad
 \begin{pmatrix}
 2 & -1 \\
 -3 & 2
 \end{pmatrix},\qquad
 \begin{pmatrix}
 2 & -1 \\
 -1 & 1
 \end{pmatrix},\qquad
 \begin{pmatrix}
 2 & 0 \\
 -n & 2
 \end{pmatrix},\\
 \begin{pmatrix}
 2 & 0 \\
 -n & 1
 \end{pmatrix},\qquad
 \begin{pmatrix}
 1 & 0 \\
 -n & 1
 \end{pmatrix}
\end{gather*}
up to a permutation of the indices.
We give several lemmas about impossible pairs.
Before giving lemmas, we note that
\begin{equation}\label{eq:Y2}
 n_{ij;p}^{+} = 0 \qquad \text{or} \qquad n_{ij;p}^{-}=0
\end{equation}
for any $i$, $j$, $p$.

\begin{Lemma}\label{lemma:ban list}
 It is impossible that $A_\pm(z) \in \mathbb{Y}$ has the following forms:
 \begin{enumerate}\itemsep=0pt
 \item[$(1)$]
 $A_+(1) =
 \bigl(\begin{smallmatrix}
 2 & -a \\
 * & *
 \end{smallmatrix}\bigr)$,
 $A_-(1) =
 \bigl(\begin{smallmatrix}
 2 & -b \\
 * & *
 \end{smallmatrix}\bigr)$ for odd $a$, $b$.
 \item[$(2)$]
 $A_+(1) =
 \bigl(\begin{smallmatrix}
 2 & -a \\
 * & *
 \end{smallmatrix}\bigr)$,
 $A_-(1) =
 \bigl(\begin{smallmatrix}
 1 & -b \\
 * & *
 \end{smallmatrix}\bigr)$ for odd $a$, $b$.
 \item[$(3)$]
 $A_+(1) =
 \bigl(\begin{smallmatrix}
 1 & -1 \\
 * & *
 \end{smallmatrix}\bigr)$,
 $A_-(1) =
 \bigl(\begin{smallmatrix}
 1 & -1 \\
 * & *
 \end{smallmatrix}\bigr)$.
 \item[$(4)$]
 $A_+(1) =
 \bigl(\begin{smallmatrix}
 1 & 0 \\
 * & *
 \end{smallmatrix}\bigr)$,
 $A_-(1) =
 \bigl(\begin{smallmatrix}
 1 & * \\
 * & *
 \end{smallmatrix}\bigr)$.
 \end{enumerate}
\end{Lemma}
\begin{proof}
 For (2), we can set
 \begin{align*}
 A_+(z) =
 \begin{pmatrix}
 1 + z^r & -f(z) \\
 * & *
 \end{pmatrix},\qquad
 A_-(z) =
 \begin{pmatrix}
 1 + z^r - z^a & -g(z) \\
 * & *
 \end{pmatrix}
 \end{align*}
 for some $f, g \in \mathbb{N}[z]$.
 By the symplectic property \eqref{eq:symplectic property}, we have
 \[
 z^a + z^{a - r} + f(z) g\bigl(z^{-1}\bigr) =
 z^{-a} + z^{r-a} + g(z) f\bigl(z^{-1}\bigr).
 \]
 Since $0 < a$ and $a-r<0$ by \eqref{eq:Y1},
 the sum of the coefficients of the terms in $f(z)g\bigl(z^{-1}\bigr)$ with positive exponents
 is equal to that with negative exponents.
 Since $f(1) g(1)$ (=$ab$) is odd, $f(z) g\bigl(z^{-1}\bigr)$ should contain the constant term $z^0$,
 which contradicts \eqref{eq:Y2}.
 The proof for (1) is similar.

 For (3), we can set
 \begin{align*}
 A_+(z) =
 \begin{pmatrix}
 1 + z^r - z^a & -z^b \\
 * & *
 \end{pmatrix},\qquad
 A_-(z) =
 \begin{pmatrix}
 1 + z^r - z^c & -z^d \\
 * & *
 \end{pmatrix}
 \end{align*}
 with $0< a,b,c,d < r$.
 Without loss of generality, we can assume $a<c$.
 By \eqref{eq:symplectic property}, we have
 \begin{equation*}
 z^{-c} + z^{r-c} + z^a + z^{a-r} + z^{c-a} + z^{d-b}
 = z^{c} + z^{c-r} + z^{-a} + z^{r-a} + z^{a-c} + z^{b-d}.
 \end{equation*}
 Since $c - a > 0$, we see that $c - a$ is equal to $c$, $r-a$, or $b-d$.
 However, the first two cases are impossible by \eqref{eq:Y1}.
 Thus $c - a = b - d$, which implies that
 \begin{equation*}
 z^{-c} + z^{r-c} + z^a + z^{a-r}
 = z^{c} + z^{c-r} + z^{-a} + z^{r-a}.
 \end{equation*}
 Since $a > 0$, we see that $a$ is equal to $c$ or $r-a$.
 However, $a = c$ is impossible by \eqref{eq:Y2}.
 Thus $a = r-a$, which implies that
$
 z^{-c} + z^{r-c}
 = z^{c} + z^{c-r}$.
 Since $c > 0$, we see that $c = r -c$.
 However, this implies that $a = r/2 = c$, which is impossible by \eqref{eq:Y2}.

 For (4), we can set
 \begin{align*}
 A_+(z) =
 \begin{pmatrix}
 1 + z^r -z^a & 0 \\
 * & *
 \end{pmatrix},\qquad
 A_-(z) =
 \begin{pmatrix}
 1 + z^r - z^b & * \\
 * & *
 \end{pmatrix}.
 \end{align*}
 By \eqref{eq:symplectic property}, we have
 \begin{align*}
 z^{a} + z^{a-r} + z^{-b} + z^{r-b} + z^{b-a} =
 z^{-a} + z^{r-a} + z^{b} + z^{b-r} + z^{a-b}.
 \end{align*}
 Comparing the number of the terms with positive and negative exponents,
 we should have $a = b$. This is impossible by \eqref{eq:Y2}.
\end{proof}

\begin{Lemma}\label{lemma:ban list indec}
 It is impossible that indecomposable $A_\pm(z) \in \mathbb{Y}$ has the form
$
 A_+(1) =
 \bigl(\begin{smallmatrix}
 * & 0 \\
 * & *
\end{smallmatrix}\bigr)$, $
 A_-(1) =
 \bigl(\begin{smallmatrix}
 * & 0 \\
 * & *
\end{smallmatrix}\bigr)$.

\end{Lemma}
\begin{proof}
 We can set
 \begin{align*}
 A_\pm(z) =
 \begin{pmatrix}
 1 + z^{r_1} - f_\pm(z) & 0 \\
 -g_\pm(z) & 1+z^{r_2} -h_\pm(z)
 \end{pmatrix}.
 \end{align*}
 Since $g_+(z)\neq 0$ or $g_-(z)\neq 0$,
 we can pick the least integer $c$ among the exponents in $g_+(z)$ and~$g_-(z)$.
 Without loss of generality, we can assume $g_+(z)$ contains the term $z^c$.
 By \eqref{eq:symplectic property}, we~have
 \begin{align}\label{eq:ban list indec}
 f_+(z) g_-\bigl(z^{-1}\bigr) + (1+z^{r_1}) g_+ \bigl(z^{-1}\bigr)=
 f_-(z) g_+\bigl(z^{-1}\bigr) + (1+z^{r_1}) g_- \bigl(z^{-1}\bigr).
 \end{align}
 The left-hand side in \eqref{eq:ban list indec} contains the term $z^{r_1 - c}$,
 but any exponent in the right-hand side is strictly smaller that $r_1 - c$ by
 \eqref{eq:Y1} and \eqref{eq:Y2}, which is a contradiction.
\end{proof}

We now search for possible pairs $A_\pm(1)$ case by case using the symplectic property
\eqref{eq:symplectic property} at~${z = 1}$ together with Lemmas~\ref{lemma:ban list} and~\ref{lemma:ban list indec}:
\begin{itemize}\itemsep=0pt
\item Case:
 $A_+(1) =
 \bigl(\begin{smallmatrix}
 2 & -1 \\
 -1 & 2
\end{smallmatrix}\bigr)$.
 The possibilities for $A_-(1)$ are
$\bigl(\begin{smallmatrix}
 2 & 0 \\
 0 & 2
\end{smallmatrix}\bigr)$,
$\bigl(\begin{smallmatrix}
 1 & 0 \\
 0 & 1
\end{smallmatrix}\bigr)$.

\item Case:
$A_+(1) =
\bigl(\begin{smallmatrix}
 2 & -1 \\
 -1 & 2
\end{smallmatrix}\bigr)$.
 The possibilities for $A_-(1)$ are $\bigl(\begin{smallmatrix}
 2 & 0 \\
 0 & 2
\end{smallmatrix}\bigr)$, $\bigl(\begin{smallmatrix}
 1 & 0 \\
 0 & 1
\end{smallmatrix}\bigr)$.

\item Case:
$A_+(1) =
\bigl(\begin{smallmatrix}
 2 & -1 \\
 -2 & 2
\end{smallmatrix}\bigr)$.
 The possibilities for $A_-(1)$ are
$\bigl(\begin{smallmatrix}
 2 & 0 \\
 -1 & 2
\end{smallmatrix}\bigr)$, $\bigl(\begin{smallmatrix}
 1 & 0 \\
 0 & 2
\end{smallmatrix}\bigr)$.

\item Case: $A_+(1) =
\bigl(\begin{smallmatrix}
 2 & -1 \\
 -3 & 2
\end{smallmatrix}\bigr)$.
 The possibilities for $A_-(1)$ are $\bigl(\begin{smallmatrix}
 2 & 0 \\
 -2 & 2
\end{smallmatrix}\bigr)$.

\item Case: $A_+(1) =
 \bigl(\begin{smallmatrix}
 2 & -1 \\
 -1 & 1
\end{smallmatrix}\bigr)$.
 The possibilities for $A_-(1)$ are $\bigl(\begin{smallmatrix}
 2 & 0 \\
 0 & 2
\end{smallmatrix}\bigr)$.

\item Case: $A_+(1) =
 \bigl(\!\begin{smallmatrix}
 2 & 0 \\
 -n & 2
\end{smallmatrix}\bigr)$.
 The possibilities for $A_\pm(1)$ are
$
\bigl(\bigl(\begin{smallmatrix}
 2 & 0 \\
 0 & 2
\end{smallmatrix}\bigr),
\bigl(\!\begin{smallmatrix}
 2 & -1 \\
 -1 & 2
 \end{smallmatrix}\bigr)\bigr)$, $
\bigl(\bigl(\begin{smallmatrix}
 2 & 0 \\
 -1 & 2
\end{smallmatrix}\bigr),
\bigl(\!\begin{smallmatrix}
 2 & -1 \\
 -2 & 2
 \end{smallmatrix}\bigr)\bigr)$, $
\bigl(\bigl(\begin{smallmatrix}
 2 & 0 \\
 0 & 2
\end{smallmatrix}\bigr),
\bigl(\begin{smallmatrix}
 1 & -1 \\
 -1 & 2
 \end{smallmatrix}\bigr)\bigr)$.

\item Case: $A_+(1) =
\bigl(\begin{smallmatrix}
 2 & 0 \\
 -n & 1
\end{smallmatrix}\bigr)$.
 The possibilities for $A_\pm(1)$ are
$\bigl( \bigl(\begin{smallmatrix}
 2 & 0 \\
 0 & 1
 \end{smallmatrix}\bigr),
\bigl(\begin{smallmatrix}
 2 & -2 \\
 -1 & 2
 \end{smallmatrix}\bigr)\bigr)$.

\item Case: $A_+(1) =
\bigl(\begin{smallmatrix}
 1 & 0 \\
 -n & 1
 \end{smallmatrix}\bigr)$.
 The possibilities for $A_\pm(1)$ are
$\bigl( \bigl(\begin{smallmatrix}
 2 & 0 \\
 0 & 1
 \end{smallmatrix}\bigr),
\bigl(\begin{smallmatrix}
 2 & -2 \\
 -1 & 2
 \end{smallmatrix}\bigr)\bigr)$.
\end{itemize}

In summary, the remaining possible pairs, up to a permutation of the indices and a change of sign,
are given in the following table:
\[
 \begin{tabular}{Sl|Sl}
 $A_+ (1)$ & $A_-(1)$ \\ \hline
 $\begin{pmatrix}
 2 & -1 \\
 -1 & 2
 \end{pmatrix}$ &
 $\begin{pmatrix}
 2 & 0 \\
 0 & 2
 \end{pmatrix}$\\
 $\begin{pmatrix}
 2 & -1 \\
 -2 & 2
 \end{pmatrix}$ &
 $\begin{pmatrix}
 2 & 0 \\
 -1 & 2
 \end{pmatrix}$\\
 $\begin{pmatrix}
 2 & -1 \\
 -3 & 2
 \end{pmatrix}$ &
 $\begin{pmatrix}
 2 & 0 \\
 -2 & 2
 \end{pmatrix}$
 \end{tabular}
 \qquad
 \begin{tabular}{Sl|Sl}
 $A_+ (1)$ & $A_-(1)$ \\ \hline
 $\begin{pmatrix}
 2 & -1 \\
 -1 & 2
 \end{pmatrix}$ &
 $\begin{pmatrix}
 1 & 0 \\
 0 & 1
 \end{pmatrix}$\\
 $\begin{pmatrix}
 2 & -1 \\
 -2 & 2
 \end{pmatrix}$ &
 $\begin{pmatrix}
 1 & 0 \\
 0 & 2
 \end{pmatrix}$\\
 $\begin{pmatrix}
 2 & -1 \\
 -1 & 1
 \end{pmatrix}$ &
 $\begin{pmatrix}
 2 & 0 \\
 0 & 2
 \end{pmatrix}$
 \end{tabular}.
\]

We now start searching for possible $A_\pm(z)$.

\begin{Lemma}\label{lemma:Apm(z) 1 2 3}
 Let $n \geq 1$.
 Suppose that
$
 A_{+}(1) =
 \bigl( \begin{smallmatrix}
 2 & - 1 \\
 - n & 2
\end{smallmatrix}\bigr)$, $
 A_{-}(1) =
 \bigl( \begin{smallmatrix}
 2 & 0 \\
 - (n-1) & 2
 \end{smallmatrix}\bigr)$.
 Then
 \begin{align*}
 A_{+}(z) =
 \begin{pmatrix}
 [2]_r & - z^{a} \\
 - z^{r-a} [n]_{2r} & [2]_{(2n-1)r}
 \end{pmatrix},\qquad
 A_{-}(z) =
 \begin{pmatrix}
 [2]_r & 0 \\
 - z^{2r-a} [n-1]_{2r} & [2]_{(2n-1)r}
 \end{pmatrix}
 \end{align*}
 for some $r$, $a$, where $[n]_r$ is the \emph{$z$-integer} defined by
$
 [n]_r \coloneqq \frac{1 - z^{rn}}{1 - z^r}$.
\end{Lemma}
\begin{proof}
 We can set
 \begin{align*}
 A_{+}(z) =
 \begin{pmatrix}
 [2]_{r_1} & - z^a \\
 - \sum\limits_{i=1}^n z^{b_i} & [2]_{r_2}
 \end{pmatrix},\qquad
 A_{-}(z) =
 \begin{pmatrix}
 [2]_{r_1} & 0 \\
 - \sum\limits_{i=1}^{n-1} z^{c_i} & [2]_{r_2}
 \end{pmatrix}.
 \end{align*}
 Without loss of generality, we can assume that
$b_1 \leq b_2 \leq \dots \leq b_n$, $
 c_1 \leq c_2 \leq \dots \leq c_{n-1}$.
 By the symplectic property \eqref{eq:symplectic property}, we have
 \begin{align*}
 \sum_{i=1}^{n-1}\bigl(z^{-c_i} + z^{r_1 - c_i}\bigr) + z^{a} + z^{a-r_2} =
 \sum_{i=1}^{n} \bigl(z^{-b_i} + z^{r_1 - b_i}\bigr).
 \end{align*}
 Comparing the degree by using the conditions \eqref{eq:Y1} and \eqref{eq:Y2}, we obtain the
 system of linear equations
$
 a = r_1 - b_1$, $ a-r_2 = -b_n$, $
 r_1 = c_i - b_i = b_{i+1} - c_i $ ($i=1,\dots,n-1$),
 which implies that
$
 r_2 = (2n-1) r_1$, $
 b_i = (2i - 1) r_1 - a$, $
 c_i = 2i r_1 - a$.
\end{proof}

\begin{Lemma}\label{lemma:Apm(z) 4}
 Suppose that
$
 A_{+}(1) =
 \bigl( \begin{smallmatrix}
 2 & - 1 \\
 -1 & 2
\end{smallmatrix}\bigr)$, $
 A_{-}(1) =
 \bigl( \begin{smallmatrix}
 1 & 0 \\
 0 & 1
\end{smallmatrix}\bigr)$.
 Then
 \begin{align*}
 A_{+}(z) =
 \begin{pmatrix}
 1+z^{2r} & - z^a \\
 -z^{2r-a} & 1+z^{2r}
 \end{pmatrix},\qquad
 A_{-}(z) =
 \begin{pmatrix}
 1+z^{2r}-z^r & 0 \\
 0 & 1+z^{2r}-z^r
 \end{pmatrix}
 \end{align*}
 for some $r$, $a$.
\end{Lemma}
\begin{proof}
 We can set
 \begin{align*}
 A_{+}(z) =
 \begin{pmatrix}
 1+z^{r_1} & - z^a \\
 -z^b & 1+z^{r_2}
 \end{pmatrix},\qquad
 A_{-}(z) =
 \begin{pmatrix}
 1+z^{r_1}-z^c & 0 \\
 0 & 1+z^{r_2}-z^d
 \end{pmatrix}.
 \end{align*}
 By \eqref{eq:symplectic property},
 we have $r_1=r_2=a+b=2c=2d$.
\end{proof}

\begin{Lemma}\label{lemma:Apm(z) 5}
 Suppose that
$
 A_{+}(1) =
 \bigl( \begin{smallmatrix}
 2 & - 1 \\
 -2 & 2
 \end{smallmatrix}\bigr)$, $
 A_{-}(1) =
 \bigl( \begin{smallmatrix}
 1 & 0 \\
 0 & 2
 \end{smallmatrix}\bigr)$.
 Then
 \begin{align*}
 A_{+}(z) =
 \begin{pmatrix}
 1+z^{2r} & - z^a \\
 -z^{2r-a}-z^{3r-a} & 1+z^{3r}
 \end{pmatrix},\qquad
 A_{-}(z) =
 \begin{pmatrix}
 1+z^{2r}-z^r & 0 \\
 0 & 1+z^{2r}
 \end{pmatrix}
 \end{align*}
 for some $r$, $a$.
\end{Lemma}
\begin{proof}
 We can set
 \begin{align*}
 A_{+}(z) =
 \begin{pmatrix}
 1+z^{r_1} & - z^a \\
 -z^{b_1}-z^{b_2} & 1+z^{r_2}
 \end{pmatrix},\qquad
 A_{-}(z) =
 \begin{pmatrix}
 1+z^{r_1}-z^c & 0 \\
 0 & 1+z^{r_2}
 \end{pmatrix}.
 \end{align*}
 Without loss of generality, we can assume $b_1 \leq b_2$.
 By \eqref{eq:symplectic property},
 we have $r_1=2c$, $r_2=3c$, $b_1 = 2c-a$, and $b_2=3c-a$.
\end{proof}

\begin{Lemma}\label{lemma:Apm(z) 6}
 Suppose that
$
 A_{+}(1) =
 \bigl( \begin{smallmatrix}
 2 & - 1 \\
 -1 & 1
 \end{smallmatrix}\bigr)$, $
 A_{-}(1) =
 \bigl( \begin{smallmatrix}
 2 & 0 \\
 0 & 2
 \end{smallmatrix}\bigr)$.
 Then
 \begin{align*}
 A_{+}(z) =
 \begin{pmatrix}
 1+z^{2r} & - z^a \\
 -z^{2r-a} & 1+z^{2r}-z^r
 \end{pmatrix},\qquad
 A_{-}(z) =
 \begin{pmatrix}
 1+z^{2r} & 0 \\
 0 & 1+z^{2r}
 \end{pmatrix}
 \end{align*}
 for some $r$, $a$.
\end{Lemma}
\begin{proof}
 We can set
 \begin{align*}
 A_{+}(z) =
 \begin{pmatrix}
 1+z^{r_1} & - z^a \\
 -z^{b} & 1+z^{r_2}-z^c
 \end{pmatrix},\qquad
 A_{-}(z) =
 \begin{pmatrix}
 1+z^{r_1} & 0 \\
 0 & 1+z^{r_2}
 \end{pmatrix}.
 \end{align*}
 By \eqref{eq:symplectic property},
 we have $r_1=r_2=a+b=2c$.
\end{proof}

\begin{proof}[Proof of Theorem \ref{theorem:main}\,(2)]
The remaining possibilities for
finite type $A_\pm (z) \in \mathbb{Y}$,
up to a~permutation of indices and change of sign,
are the six families of the pairs
given in Lemmas~\ref{lemma:Apm(z) 1 2 3}--\ref{lemma:Apm(z) 6}, which contain the parameters $r$, $a$.
We can verify that these six families belong to $\mathbb{Y}$
if only if
\begin{itemize}\itemsep=0pt
 \item $1 < r$ and $0 < a < r$ for Lemma \ref{lemma:Apm(z) 1 2 3},
 \item $1 < r$ and $0 < a < 2r$ for Lemmas~\ref{lemma:Apm(z) 4}, \ref{lemma:Apm(z) 5}, and \ref{lemma:Apm(z) 6},
\end{itemize}
where the only if part follows from \eqref{eq:Y1}.
With these conditions,
the families for $n = 1, 2, 3$ in Lemma~\ref{lemma:Apm(z) 1 2 3}
can be reduced to the pairs (1), (2), and (3) in Table \ref{table:finie type}, respectively,
by change of slices,
and the families in Lemmas~\ref{lemma:Apm(z) 4}, \ref{lemma:Apm(z) 5}, and \ref{lemma:Apm(z) 6}
can be reduced to the pairs (4), (5), and (6) in Table \ref{table:finie type}, respectively,
by change of slices.
\end{proof}

We note that the pairs in Table~\ref{table:finie type} correspond to the
value $(r, a) = (2, 1)$ in Lemmas~\ref{lemma:Apm(z) 1 2 3}--\ref{lemma:Apm(z) 6}.

\subsection*{Acknowledgements}
 The authors would like to thank the anonymous referees for their constructive and helpful suggestions.
	This work is supported by JSPS KAKENHI Grant Number JP21J00050.

\pdfbookmark[1]{References}{ref}
\LastPageEnding

\end{document}